 
\input amstex
\input epsf
\loadmsbm
\documentstyle{gen-j}
\NoBlackBoxes
\font\tit=cmr10 scaled\magstep3
\font\sfa=cmss10 scaled\magstep1

\def\a{\vrule height 11.25pt  depth 4.75pt} 
\def\({[}
\def\){]}
 \catcode`\@=11
\def\bqed{\ifhmode\unskip\nobreak\fi\quad
  \ifmmode\blacksquare\else$\m@th\blacksquare$\fi}
\def\mex{\qopname@{mex}}

\catcode`\@=\active
\UseAMSsymbols
\vglue1.0cm
\topmatter
\title{}\endtitle

\endtopmatter

\centerline{{\tit Heap Games, Numeration Systems and Sequences}}
\vskip1.2cm
\centerline{\sfa Aviezri S.\ Fraenkel}\medskip
\centerline{Department of Applied Mathematics and Computer Science} 
\centerline{Weizmann Institute of Science} 
\centerline{Rehovot 76100, Israel}
\centerline{\tt fraenkel\@wisdom.weizmann.ac.il}
\centerline{\tt http://www.wisdom.weizmann.ac.il/\~\! fraenkel}

\vskip1.0cm
{\bf Abstract.} We propose and analyse a 2-parameter family of 2-player 
games on two heaps of tokens, and present a strategy based on a class 
of sequences. The strategy looks easy, but is actually hard. A class of 
exotic numeration systems is then used, which enables us to decide 
whether the family has an efficient strategy or not. We introduce 
yet another class of sequences, and demonstrate its equivalence 
with the class of sequences defined for the strategy of our games.

\vskip0.8cm
\document
\centerline{\bf 1. An Example}\medskip

Given a two-player game played on two heaps (piles) of finitely many tokens. 
There are two types of moves: I. Take any positive number of tokens 
from {\sl one\/} heap, possibly the entire heap. II. Take from {\sl both\/} 
heaps, $k$ from one and $l$ from the other, with, say, $k\le l$. Then 
the move is constrained by the condition $0<k\le l<2k+2$, which is 
equivalent to $0\le l-k<k+2,\ k>0$. The player making the last move 
(after which both heaps are empty) wins, and the opponent loses. 

A position $q$ in a game of this sort is called a $P$-position, if the 
{\it P}revious player can win, i.e., the player who moved to $q$. 
It is an $N$-position, if the {\it N}ext player can win, i.e., the 
player moving from $q$. The position $(0,0)$ (two empty heaps) is a 
$P$-position, since the first player cannot even make a move, so the second 
wins by default. The next $P$-position is $(1,4)$: if Jean takes an 
entire heap, then Gill takes the other and wins. If Jean takes any 
part of the larger heap, Gill can take the balance of both heaps. 
Lastly, Jean cannot remove both heaps, and if she takes from both 
heaps, then Gill takes the balance and wins. 

Table~1 lists the first few $P$-positions. The reader will do well to 
try and construct the next few entries of the table before reading on. 

\midinsert
\topcaption{Table 1. The first few $P$-positions.}
\endcaption
$$\vcenter{\lineskip -2pt\halign{\tabskip 10pt\hfill #\hfill\ &\hfill #\ &
\hfill #\hfill\ &\hfill #\ \ &\hfill #\ \ &\tabskip 0pt\hfill #\hfill\cr
\noalign{\hrule}
\a&$n$&\a&$A_n$&$B_n$&\a\cr
\noalign{\hrule}
\a&0&\a&0&0&\a\cr
\a&1&\a&1&4&\a\cr
\a&2&\a&2&8&\a\cr
\a&3&\a&3&12&\a\cr
\a&4&\a&5&18&\a\cr
\a&5&\a&6&22&\a\cr
\a&6&\a&7&26&\a\cr
\a&7&\a&9&32&\a\cr
\a&8&\a&10&36&\a\cr
\a&9&\a&11&40&\a\cr
\a&10&\a&13&46&\a\cr
\a&11&\a&14&50&\a\cr
\a&12&\a&15&54&\a\cr
\a&13&\a&16&58&\a\cr
\noalign{\hrule}}}$$
\endinsert

If $S$ is any finite subset of nonnegative integers, denote by 
$\mex S$ the least nonnegative integer in the complement of $S$, 
i.e., the least nonnegative integer not occurring in $S$. Note 
that the $\mex$ of the empty set is $0$. The term $\mex$, used in 
\cite{BCG1982}, stands for Minimum EXcluded value. The structure 
of Table~1 is made explicit by:

$$A_n=\mex\{A_i,B_i:i<n\},\ B_n=2(A_n+n)\quad (n\ge 0).$$

This is a special case of Theorem~1 below, in the proof of which we also 
see that if $A=\bigcup_{n=1}^{\infty} A_n$, $B=\bigcup_{n=1}^{\infty} B_n$, 
then $A$ and $B$ are {\it complementary\/}, i.e., $A\cup B=$ set of all 
positive integers, and $A\cap B=\emptyset$. 

Given any two heaps of our game, containing $x$ and $y$ tokens with 
$x\le y$. The complementarity of $A$ and $B$ implies that either $x=A_n$ 
or $x=B_n$ for some $n$. Hence Table~1 has to be computed only up to 
the encounter of $x$. Moreover, it is not hard to see that $n\le x$, 
and if $x=A_n$, then $x/2<n$, so the table has to be computed up to at 
most $\Omega (x)$, which implies a strategy computation linear in $x$, 
which looks good. 

The trouble with this strategy is the same as that of the simple-minded 
primality-testing algorithm for a given integer $m$: divide $m$ by the 
integers $\le\sqrt{m}$, and if none of them divides $m$, then $m$ is prime. 
This algorithm is linear in $m$. The problem in both cases is of course 
that the input length is the $\log$ of the input numbers $x$, $y$ and $m$, 
rather than $x$, $y$ and $m$ themselves. 

The two algorithms mentioned above, that for the strategy computation, 
and that for primality testing are thus actually exponential in the 
input length. 

The central question we address here is whether games of the type 
considered above have a polynomial strategy, or whether their best 
strategies are necessarily exponential. Before that we define in 
\S2 the family of games precisely, introduce a family of sequences, 
and formulate and prove the winning strategy in terms of these sequences.

In \S3 we present an argument against polynomiality of the games, and 
in \S4 we introduce a numeration system that turns out to be relevant 
to our games. The connection between the games and the numeration system 
is made explicit in \S5. This enables us to decide the games' polynomiality 
question in \S6. Yet another class of sequences is introduced in \S7, 
where we prove equivalence between the two classes of sequences. In 
the final \S8 we summarize our results, give motivation and present 
a few open problems.\bigskip

\centerline{\bf 2. A Family of Heap Games and their Winning Strategies}\medskip

Denote by $\Bbb Z^0$ and $\Bbb Z^+$ the set of nonnegative integers and 
positive integers respectively. Our family of heap games depends on 
two parameters $s, t\in\Bbb Z^+$. Given are two heaps of 
finitely many tokens. There are two types of moves: I. Take any 
positive number of tokens from a single heap, possibly the entire 
heap. II. Take $k>0$ and $l>0$ from the two heaps, say $0<k\le l$. 
This move is constrained by the condition

$$0<k\le l<sk+t,\tag 1$$
which is equivalent to $0\le l-k<(s-1)k+l$,\quad $k\in\Bbb Z^+$. 

The example presented in \S1 is the special case $s=t=2$. Denote by 
$\Cal P$ the set of all $P$-positions. 

\proclaim{\bf Theorem 1} $\Cal P=\bigcup_{i=0}^{\infty}\{(A_i,B_i)\}$, 
where\/ 
$$A_n=\mex\{A_i, B_i:0\le i<n\},\ B_n=sA_n+tn\quad (n\in\Bbb Z^0).\tag 2$$
\endproclaim 

\demo{\bf Proof} Let $A=\bigcup_{n=1}^{\infty} A_n$, 
$B=\bigcup_{n=1}^{\infty} B_n$. 
Then $A$, $B$ are complementary with respect to $\Bbb Z^+$:\ 
$A\cup B=\Bbb Z^+$ follows from the $\mex$ property. Suppose that 
$A_m=B_n$. Then $m>n$ implies that $A_m$ is the $\mex$ of a set 
containing $B_n=A_m$, a contradiction. If $m\le n$, then 
$B_n=sA_n+tn\ge sA_m+tm>A_m$, another contradiction. Thus $A\cap B=\emptyset$.
We will also need the fact that $A_n$ and $B_n$ are strictly increasing 
sequences, which is clear from their definition. 

Let $W=\bigcup_{i=0}^{\infty} (A_i,B_i)$. It evidently suffices to show 
two things: I.~~A player moving from some $(A_n,B_n)\in W$ lands in a 
position not in $W$. II.~~Given any position $(x,y)$ not in $W$, there 
is a move to some $(A_n,B_n)\in W$.\medskip 

I.~~A move of the first type from $(A_n,B_n)\in W$ clearly leads to a
position not in $W$, since $A_n$ and $B_n$ are strictly increasing, so 
they have no repeating terms. Suppose that a move of the second type 
from $(A_n,B_n)\in W$ produces a position $(A_m,B_m)\in W$. 
Then $m<n$. For $k=A_n-A_m$, $l=B_n-B_m$, we have 
$$l=sA_n+tn-sA_m-tm=s(A_n-A_m)+t(n-m)\ge sk+t,$$ 
which contradicts condition (1).\medskip    

II.~~Let $(x,y)$ with $x\le y$ be a position not in $W$. Since $A$ and 
$B$ are complementary, every positive integer appears exactly once in 
exactly one of $A$ and $B$. Therefore we have either $x=B_n$ or else 
$x=A_n$ for some $n\ge 0$.\medskip 

{\it Case\/} (i): $x=B_n$. Then move $y\rightarrow A_n$.

{\it Case\/} (ii): $x=A_n$. If $y>B_n$, then move $y\rightarrow B_n$. 
So suppose that $A_n\le y<B_n$. If $y<sA_n+t$, move $(x,y)\rightarrow (0,0)$, 
which satisfies (1) with $k=A_n$, $l=y$. So let $y\ge sA_n+t$. Put 
$m=\lfloor (y-sA_n)/t\rfloor$, and move $(x,y)\rightarrow (A_m,B_m)$, 
where $\lfloor x\rfloor$ denotes the largest integer $\le x$. This 
move is legal, since (a) $m<n$,\ \ (b) $y>B_m$,\ \ 
(c) $A_n-A_m\le y-B_n<s(A_n-A_m)+t$. Indeed,\medskip   

(a) $y-sA_n<B_n-sA_n=tn$, so 
$m=\lfloor (y-sA_n)/t\rfloor\le (y-sA_n)/t<n$;

(b) $m\le(y-sA_n)/t $, so $y\ge tm+sA_n=B_m+s(A_n-A_m)>B_m$;

(c) $m>((y-sA_n)/t)-1$, so $y<tm+sA_n+t$; by (b), $y-B_m\ge A_n-A_m$, 
hence 
$$A_n-A_m\le y-B_m<tm+sA_n+t-sA_m-tm=s(A_n-A_m)+t,$$
and (1) is satisfied.\bqed\enddemo\medskip

The {\sl statement\/} of Theorem~1 enables one to decide whether any 
given position $(x,y)$ is a $P$-position or $N$-position, and the 
{\sl proof\/} clearly indicates a winning move from any $N$-position. 
These two things together constitute a winning strategy for the game. 

However, as was pointed out in \S1 after Table~1, the strategy is 
exponential (the inequalities for $x$ hold for all $s, t\in\Bbb Z^+$, 
not just for the special example considered there). But 
only the construction of the table needs exponential time and, in fact, 
exponential space. The rest of the algorithm is polynomial.\bigskip

\centerline{\bf 3. An Argument Against Polynomiality}\medskip

Suppose that there exist real numbers $\alpha$ and $\beta$ such 
that for the $A_n$ and $B_n$ defined in Theorem~1, 
$A_n=\lfloor n\alpha\rfloor$, and $B_n=\lfloor n\beta\rfloor$ 
for all $n\in\Bbb Z^0$. A simple density argument then shows 
that $\alpha$ and $\beta$ must satisfy $\alpha^{-1}+\beta^{-1}=1$, 
hence $1<\alpha<2<\beta$, and $\alpha$, $\beta$ are in fact irrational. 

A strategy based on this observation can be applied to any given 
game position $(x,y)$. We have 
$$\multline
x=\lfloor n\alpha\rfloor\iff x<n\alpha<x+1\\
\iff\frac x{\alpha}<n
<\frac{x+1}{\alpha}\iff\left\lfloor\frac{x+1}{\alpha}\right\rfloor=
\left\lfloor\frac x{\alpha}\right\rfloor+1.
\endmultline$$ 

Therefore either $x=\lfloor n\alpha\rfloor =A_n$ where 
$n=\lfloor (x+1)/\alpha\rfloor$, or else, by complementarity, 
$x=\lfloor n\beta\rfloor =B_n$, where $n=\lfloor (x+1)/\beta\rfloor$. 
We have thus reduced the situation to that considered in cases~(ii) and 
(i) in the proof of Theorem~1, and hence the strategy presented in that 
proof can be followed. For example, if $x=\lfloor n\alpha\rfloor =A_n$ 
and $s\lfloor n\alpha\rfloor +t\le y<s\lfloor n\alpha\rfloor+tn=
\lfloor n\beta\rfloor$, then for 
$m=\lfloor (y-s\lfloor n\alpha\rfloor)/t\rfloor$, we move 
$(x,y)\rightarrow (\lfloor m\alpha\rfloor ,\lfloor m\beta\rfloor)\in\Cal P$. 
For implementing this strategy, $\alpha$ has to be computed to a precision 
of $O(\log x)$ digits, and its storage requires $O(\log x)$ words, 
which is linear in the input size of $x$ (given in binary, say). 
Thus this strategy {\sl is\/} polynomial. See also the remark at the 
end of the previous section. 

Is it far-fetched to hope for the existence of such real numbers 
$\alpha$ and $\beta$? Well, for the special case $s=t=1$ our games 
reduce to Wythoff's game, for which such real numbers indeed 
exist, namely $\alpha=(1+\sqrt{5})/2$, $\beta=(3+\sqrt{5})/2$. This 
was already shown in \cite{Wyt07}. See also \cite{Cox53}, \cite{YaYa67}. 
In \cite{Fra82} a generalization of Wythoff's game was proposed, 
namely the case of any $t\in\Bbb Z^+$, but $s=1$. Also for this case 
these numbers exist, namely, $\alpha=(2-t+\sqrt{t^2+4})/2$, 
$\beta=\alpha +t$. 

We now show, however, that for $s>1$, such real numbers cannot exist! 

\proclaim{\bf Theorem 2} For\/ $A_n$, $B_n$ as defined in Theorem~\rm{1}, 
there exist real numbers\/ $\alpha$, $\gamma$, $\beta$, $\delta$ 
such that\/ $A_n=\lfloor n\alpha +\gamma\rfloor$ and\/ 
$B_n=\lfloor n\beta +\delta\rfloor$ for all $n\in\Bbb Z^0$, if and 
only if\/ $s=1$.\endproclaim 

\demo{\bf Proof} Since the sequence $\{A_n\}$ is strictly increasing, 
$$B_{n+1}-B_n=sA_{n+1}+t(n+1)-sA_n-tn=s(A_{n+1}-A_n)+t\ge s+t\ge 2.$$
Since $B_n=sA_n+tn$, the sequence $\{B_n\}$ is nonempty. Since $A$ and $B$ 
are complementary, we thus cannot have $A_{n+1}-A_n=1$ for all $n$. 
Therefore there exists $n$ such that $A_{n+1}-A_n\ge 2$. Hence there 
is $n$ for which $B_{n+1}-B_n\ge 2s+t\ge 3$. It follows that there is 
$n$ for which $A_{n+1}-A_n=1$. Since for all $n\in\Bbb Z^0$ we have 
$B_{n+1}-B_n\ge 2$, there can be no $n$ for which $A_{n+1}-A_n>2$. Hence 
$A_{n+1}-A_n\in\{1,2\}$ and $B_{n+1}-B_n\in\{s+t,2s+t\}$ for all $n\ge 0$. 

Given a nondecreasing sequence of integers $S=a_1,a_2,\dots$, the 
{\it spectrum\/} question is whether there exist real 
numbers $\alpha$, $\gamma$ such that $S=\lfloor n\alpha+\gamma\rfloor$.
The spectrum terminology is used in \cite{GLL78}, 
where it is shown, for the {\sl homogeneous\/} case ($\gamma =0$), that 
if the prefix $M_r$ of length $r$ of $S$ is ``nearly linear'', then it 
is the beginning of a spectrum. If it is, we'll say that $M_r$ is 
{\it spectral}. In \cite{BoFr81}, necessary and sufficient conditions 
are given for $M_r$ to be spectral in the (possibly) nonhomogeneous 
case. See also \cite{BoFr84}. 

Let 
$$\underline{d}(M_r)=\max_{1\le i<k\le r} \frac{a_k-a_{k-i}-1}i~,\quad 
\overline{d}(M_r)=\min_{1\le i<k\le r} \frac{a_k-a_{k-i}+1}i~.$$
One of the necessary and sufficient conditions for $M_r$ to be spectral 
given in \cite{BoFr81} is that $\underline{d}(M_r)<\overline{d}(M_r)$. 

Put $a_k=B_{n+1}$, $a_{k-1}=B_n$. For any portion of length $r$ of 
the sequence $\{B_n\}$ for which both the difference $2s+t$ and $s+t$ occurs, 
we then have, where we use the larger difference for $\underline{d}$ 
and the smaller for $\overline{d}$, 
$$\underline{d}(M_r)\ge a_k-a_{k-1}-1=2s+t-1,\quad 
\overline{d}(M_r)\le a_k-a_{k-1}+1=s+t+1.$$
So a necessary condition for $M_r$ to be spectral is $2s+t-1<s+t+1$, 
which holds if and only if $s<2$, i.e., $s=1$. For $s=1$ the sequence 
$\{B_n\}$ is indeed a spectrum, as remarked above. 

The structure of $\{B_n\}$ implies that in $\{A_n\}$ there are runs 
of $1$s of length $s+t-2$, and $2s+t-2$. An argument analogous to the 
above then leads to the necessary condition $2s+t-3<s+t-1$, which again 
leads to $s=1$, for which case indeed $\{A_n\}$ is a spectrum.
\bqed\enddemo

Thus the question whether our heap games have a polynomial strategy 
or not is still open for $s>1$.\bigskip

\centerline{{\bf 4. A Class of Exotic Numeration Systems}}\medskip

For $u_{-1}$ a constant, $u_0=1$ and $b_1$, $b_2$ integers satisfying 
$b_1\ge b_2\ge 1$, consider the linear recurrence 
$u_n=b_1u_{n-1}+b_2u_{n-2}\ (n\ge 1)$. We can consider the $u_0,u_1,\dots$ 
as bases of a numeration system with digits $d_i\in\{0,\dots,b_1\}$. 
But then an integer such as $u_n$ has two representations: $u_n$ 
itself, and $b_1u_{n-1}+b_2u_{n-2}$. Since we would like to have 
uniqueness of representation, it is natural to require that 
$d_i=b_1\implies d_{i-1}<b_2$\ $(i\ge 1)$. It turns out that under this 
condition every positive integer $m$ indeed has a unique representation. 
This is a special case of Theorem~2 in \cite{Fra85}. Moreover, the greedy 
algorithm of repeatedly dividing $m$ or its remainder by the largest 
$u_i$ not exceeding this remainder, yields this unique representation. 
The case $b_1=b_2=1$ gives a binary representation known as the 
Zeckendorf representation \cite{Zec72}.\medskip 

{\bf Example.} We consider the case $u_{-1}=\frac12$, $(b_1,b_2)=(3,2)$. 
Then $u_1=4$, $u_2=14$, $u_3=50$, $u_4=178$, $\dots$ . The representations 
of the integers 1 to 60 in this numeration system are displayed in 
Table~2.\medskip

\midinsert
\topcaption{Table 2. Representation of first few integers in $\Bbb Z^+$.}
\endcaption
$$\vcenter{\lineskip -2pt\halign{\tabskip 15pt\hfill #\hfill&\hfill #
\hfill&\hfill #\hfill&\hfill #\hfill&\hfill #\hfill&\hfill #\hfill&
\hfill #\hfill&\hfill #\hfill&\hfill #\hfill&\hfill #\hfill&
\hfill #\hfill&\hfill #\hfill&\hfill #\hfill&\tabskip 0pt\hfill #\hfill\cr
\noalign{\hrule}
\a&50&14&4&1&\a&$n$&\a$\,$\a&14&4&1&\a&$n$&\a\cr
\noalign{\hrule}
\a&&2&0&3&\a&31&\a$\,$\a&&&1&\a&1&\a\cr
\a&&2&1&0&\a&32&\a$\,$\a&&&2&\a&2&\a\cr
\a&&2&1&1&\a&33&\a$\,$\a&&&3&\a&3&\a\cr
\a&&2&1&2&\a&34&\a$\,$\a&&1&0&\a&4&\a\cr
\a&&2&1&3&\a&35&\a$\,$\a&&1&1&\a&5&\a\cr
\a&&2&2&0&\a&36&\a$\,$\a&&1&2&\a&6&\a\cr
\a&&2&2&1&\a&37&\a$\,$\a&&1&3&\a&7&\a\cr
\a&&2&2&2&\a&38&\a$\,$\a&&2&0&\a&8&\a\cr
\a&&2&2&3&\a&39&\a$\,$\a&&2&1&\a&9&\a\cr
\a&&2&3&0&\a&40&\a$\,$\a&&2&2&\a&10&\a\cr
\a&&2&3&1&\a&41&\a$\,$\a&&2&3&\a&11&\a\cr
\a&&3&0&0&\a&42&\a$\,$\a&&3&0&\a&12&\a\cr
\a&&3&0&1&\a&43&\a$\,$\a&&3&1&\a&13&\a\cr
\a&&3&0&2&\a&44&\a$\,$\a&1&0&0&\a&14&\a\cr
\a&&3&0&3&\a&45&\a$\,$\a&1&0&1&\a&15&\a\cr
\a&&3&1&0&\a&46&\a$\,$\a&1&0&2&\a&16&\a\cr
\a&&3&1&1&\a&47&\a$\,$\a&1&0&3&\a&17&\a\cr
\a&&3&1&2&\a&48&\a$\,$\a&1&1&0&\a&18&\a\cr
\a&&3&1&3&\a&49&\a$\,$\a&1&1&1&\a&19&\a\cr
\a&1&0&0&0&\a&50&\a$\,$\a&1&1&2&\a&20&\a\cr
\a&1&0&0&1&\a&51&\a$\,$\a&1&1&3&\a&21&\a\cr
\a&1&0&0&2&\a&52&\a$\,$\a&1&2&0&\a&22&\a\cr
\a&1&0&0&3&\a&53&\a$\,$\a&1&2&1&\a&23&\a\cr
\a&1&0&1&0&\a&54&\a$\,$\a&1&2&2&\a&24&\a\cr
\a&1&0&1&1&\a&55&\a$\,$\a&1&2&3&\a&25&\a\cr
\a&1&0&1&2&\a&56&\a$\,$\a&1&3&0&\a&26&\a\cr
\a&1&0&1&3&\a&57&\a$\,$\a&1&3&1&\a&27&\a\cr
\a&1&0&2&0&\a&58&\a$\,$\a&2&0&0&\a&28&\a\cr
\a&1&0&2&1&\a&59&\a$\,$\a&2&0&1&\a&29&\a\cr
\a&1&0&2&2&\a&60&\a$\,$\a&2&0&2&\a&30&\a\cr
\noalign{\hrule}}}$$
\endinsert

A question we just might ask at this point is whether there is any 
connection between Tables~1 and 2. If we scan the first few entries of both, 
we may be tempted to conclude that the entries under $A_n$ in Table~1 all have 
representation ending in no $0$. But then $14$ is a counterexample, 
whose representation ends in two $0$s. Also it appears that the 
$B_n$ all have representation ending in a single $0$. But $50$, 
with representation $1000$ is a counterexample, in fact, the only 
counterexample in the range of the two tables. 

However, there is no counterexample, as far as the two tables go, to 
the following two remarkable, \ae sthetically pleasing, properties:\medskip

{\bf a.} All the $A_n$ have representations ending in an {\sl even\/} 
number of $0$s, and all the $B_n$ have representations ending in an 
{\sl odd\/} number of $0$s. 

{\bf b.} For every $(A_n,B_n)\in\Cal P$, the representation of 
$B_n$ is the ``left shift'' of the representation of $A_n$.\medskip 

Thus $(1,4)$ of Table~1 has representation $(1,10)$, and $(6,22)$ 
has representation $(12,120)$: $10$ is the ``left shift'' of $1$, 
$120$ the left shift of $12$. We remark that the second part of 
{\bf a} is not independent; it follows from its first part, since 
$A$ and $B$ are complementary. 

In the next section we state these properties in a precise manner 
and prove their validity.\bigskip

\centerline{{\bf 5. Wedding Numeration Systems with Heap Games}}\medskip

For fixed $s, t\in\Bbb Z^+$, put $u_{-1}=1/s$, $u_0=1$, and let 
$u_n=(s+t-1)u_{n-1}+su_{n-2}\ (n\ge 1)$. Denote by $\Cal U$ the numeration 
system with bases $u_0$, $u_1$, $\ldots$ and digits $d_i\in\{0,\dots,s+t-1\}$ 
such that $d_{i+1}=s+t-1\implies d_i<s\ (i\ge 0)$. Every positive 
integer has a unique representation over $\Cal U$, as mentioned in 
the previous section.\medskip 

{\bf Notation and Definitions.}
\roster
\item"{(a)}" For every $m\in\Bbb Z^0$ write $R(m)$ for the representation 
of $m$ over $\Cal U$. 
\item"{(b)}" Denote by $LR(m)$ the ``left shift'' of $R(m)$, i.e., if 
$R(m)=\sum_{i=0}^{n} d_iu_i$, then $LR(m)=\sum_{i=0}^{n} d_iu_{i+1}$. 
\item"{(c)}" A positive integer $m$ is {\it EVen-taILed\/} (for short:
{\it evil\/}), if $R(m)$ ends in an even (possibly $0$) number of $0$s. 
It is {\it Odd-taiLeD\/} (for short: {\it old\/}), if $R(m)$ ends in an 
odd number of $0$s. It is convenient to let $0$ be both evil and old. 
\item"{(d)}" Put $q=s-1$, and $r=s+t-1$. Then the above recurrence has 
the form $u_{-1}=1/s$, $u_0=1$, $u_n=ru_{n-1}+su_{n-2}\ (n\ge 1)$; 
and the representation with digits $d_i\in\{0,\dots,r\}$ satisfies 
$d_{i+1}=r\implies d_i\le q\ (i\ge 0)$.
\endroster

We mention that in \cite{BCG82, \rm Ch.~4}, ``evil number'' 
is used for a number whose binary expansion contains an even number of 
$1$s ({\it even weight\/} in coding theory language). 

\proclaim{\bf Lemma 1} For\/ $m\in\Bbb Z^+$, let\/ 
$R(m)=\sum_{i=0}^{n} d_iu_i$. 

\rm{(i)} Suppose that for some\/ $k\in\Bbb Z^0$, the tail of\/ $R(m)$ 
has digits 
$$d_{2k}d_{2k-1}d_{2k-2}\dots d_3d_2d_1d_0=d_{2k}rq\dots rqrq,$$
where\/ $d_{2k}\in\{0,\dots,q\}$ and\/ $d_{2k}=q\implies d_{2k+1}<r$. 
Then\/ $R(m+1)=(d_{2k}+1)u_{2k}+\sum_{i=2k+1}^{n} d_iu_i$, so $m+1$ is evil. 

\rm{(ii)} Suppose that for some\/ $k\in\Bbb Z^0$, the tail of\/ $R(m)$ 
has digits 
$$d_{2k+1}d_{2k}d_{2k-1}\dots d_2d_1d_0=d_{2k+1}rq\dots rqr,$$
where\/ $d_{2k+1}\in\{0,\dots,q\}$ and $d_{2k+1}=q\implies d_{2k+2}<r$. 
Then\/ $R(m+1)=(d_{2k+1}+1)u_{2k+1}+\sum_{i=2k+2}^{n} d_iu_i$, so $m+1$ 
is old. 
\endproclaim

{\bf Note.} We point out the special case $k=0$, where for (i), 
$d_0\le q$ and $d_0=q\implies d_1<r$; for (ii), $d_1\le q$ and 
$d_1=q\implies d_2<r$. 

\demo{\bf Proof} We note that the hypothesis on $d_{2k}$ for (i) implies 
that $(d_{2k}+1)u_{2k}+\sum_{i=2k+1}^{n} d_iu_i$ is a legal representation 
over $\Cal U$. Similarly for (ii).\medskip 

(i) By adding and subtracting $u_0$, we get 
$$\split
m&=(qu_0+ru_1)+(qu_2+ru_3)+\dots +(qu_{2k-2}+ru_{2k-1})+d_{2k}u_{2k}\\
&+\sum_{i=2k+1}^{n} d_iu_i=(d_{2k}+1)u_{2k}+\sum_{i=2k+1}^{n} d_iu_i -1.
\endsplit$$
Thus $m+1=(d_{2k}+1)u_{2k}+\sum_{i=2k+1}^{n} d_iu_i=R(m+1)$. 

(ii) Adding and subtracting $su_{-1}=1$, gives 
$$\split
m&=ru_0+(qu_1+ru_2)+(qu_3+ru_4)+\dots +(qu_{2k-1}+ru_{2k})\\
&+d_{2k+1}u_{2k+1}
+\sum_{i=2k+2}^{n} d_iu_i=(d_{2k+1}+1)u_{2k+1}+\sum_{i=2k+2}^{n} d_iu_i -1.
\endsplit$$
Thus $m+1=(d_{2k+1}+1)u_{2k+1}+\sum_{i=2k+2}^{n} d_iu_i=R(m+1)$.\bqed
\enddemo\medskip

\proclaim{Lemma 2} Consider the set of pairs\/ 
$\bigcup_{k=0}^{\infty} (V_k,W_k)$, where\/ $0=V_0<V_1<\dots$ 
is the set of all evil numbers, and\/ $R(W_k)=LR(V_k)$ for all\/ $k$. 
Then\/ $W_k-sV_k=tk$ for all\/ $k$.\endproclaim 

\demo{\bf Proof} Induction on $k$. The assertion holds trivially for 
$k=0$. Suppose that $W_k-sV_k=tk$ for arbitrary fixed $k$. Let 
$R(V_k)=\sum_{i=0}^{n} d_iu_i$. Then 
$$R(W_k)-sR(V_k)=LR(V_k)-sR(V_k)=\sum_{i=0}^{n} d_i(u_{i+1}-su_i),$$
so 
$$tk=\sum_{i=0}^{n} d_i(u_{i+1}-su_i).\tag 3$$
We consider three cases. 

{\bf I}.~~The tail of $R(V_k)$ is as in case~(i) of Lemma~1. Then $V_k+1$ is 
evil, so $V_{k+1}=V_k+1$ and 
$R(V_{k+1})=(d_{2k}+1)u_{2k}+\sum_{i=2k+1}^{n} d_iu_i$. Thus 
$$\multline 
LR(V_{k+1})-sR(V_{k+1})=(d_{2k}+1)(u_{2k+1}-su_{2k})\\
+\sum_{i=2k+1}^{n} d_i(u_{i+1}-su_i)
=u_{2k+1}-su_{2k}+\sum_{i=2k}^{n} d_i(u_{i+1}-su_i).
\endmultline\tag 4$$

For case~(i) of Lemma~1 we have by (3), 
$$\split 
tk&=q(u_1-su_0)+r(u_2-su_1)+q(u_3-su_2)+r(u_4-su_3)\\
&+\dots +q(u_{2k-1}-su_{2k-2})+r(u_{2k}-su_{2k-1})
+\sum_{i=2k}^{n} d_i(u_{i+1}-su_i). 
\endsplit$$
We sum together the positive terms, adding and subtracting $u_1$. Then 
$ru_2+su_1=u_3$ is added to $(s-1)u_3$, and so on, leading 
to $u_{2k+1}-u_1$. We then sum all the negative terms, subtracting 
and adding $su_0$, leading to $-su_{2k}+su_0$. Thus 
$$\split
tk&=u_{2k+1}-u_1-su_{2k}+su_0+\sum_{i=2k}^{n} d_i(u_{i+1}-su_i)\\
&=u_{2k+1}-su_{2k}-t+\sum_{i=2k}^{n} d_i(u_{i+1}-su_i).
\endsplit$$
Hence by (4), 
$$t(k+1)=u_{2k+1}-su_{2k}+\sum_{i=2k}^{n} d_i(u_{i+1}-su_i)=
LR(V_{k+1})-sR(V_{k+1}),$$ 
as was to be shown.\medskip

{\bf II}.~~The tail of $R(V_k)$ is as in case~(ii) of Lemma~1. Then $V_k+1$
is old, but $V_k+2$ is clearly evil, since $R(V_k+2)$ ends in $1$, so 
$V_{k+1}=V_k+2$. Then $V_{k+1}=1+(d_{2k+1}+1)u_{2k+1}
+\sum_{i=2k+2}^{n} d_iu_i$, 
so $R(V_{k+1})=u_0+(d_{2k+1}+1)u_{2k+1}+\sum_{i=2k+2}^{n} d_iu_i$. Hence 
$$\multline
LR(V_{k+1})-sR(V_{k+1})=(u_1-su_0)+(d_{2k+1}+1)(u_{2k+2}-su_{2k+1})\\
+\sum_{i=2k+2}^{n} d_i(u_{i+1}-su_i)
=t+u_{2k+2}-su_{2k+1}+\sum_{i=2k+1}^{n} d_i(u_{i+1}-su_i).
\endmultline\tag 5$$
For case~(ii) of Lemma~1 we have by (3), 
$$\split 
tk&=r(u_1-su_0)+q(u_2-su_1)+r(u_3-su_2)+q(u_4-su_3)\\
&+\dots +q(u_{2k}-su_{2k-1})+r(u_{2k+1}-su_{2k})
+\sum_{i=2k+1}^{n} d_i(u_{i+1}-su_i). 
\endsplit$$
Summing the positive terms, adding and subtracting $su_0$, 
leads to $u_{2k+2}-s$. Summing the negative terms, subtracting and 
adding $s^2q_{-1}=s$, gives $-su_{2k+1}+s$. Thus $tk=u_{2k+2}-su_{2k+1}+
\sum_{i=2k+1}^{n} d_i(u_{i+1}-su_i)$.

Thus by (5), $LR(V_{k+1})-sR(V_{k+1})=t(k+1)$, as required.\medskip

{\bf III}.~~The digit $d_0$ satisfies $q<d_0<r$. Since clearly $d_1<r$, 
we have $V_{k+1}=V_k+1$, so by (3), 
$$\split
LR(V_{k+1})-sR(V_{k+1})&=(d_0+1)(u_1-su_0)+\sum_{i=1}^{n} d_i(u_{i+1}-su_i)\\
&=t+\sum_{i=0}^{n} d_i(u_{i+1}-su_i)=t(k+1).
\endsplit$$
Every evil number $P$ is of one of the three forms considered above. Note 
that if $P$ ends in an even number of $0$s, it is of the form {\bf I}. 
\bqed\enddemo\medskip

Lemma 2 enables us to prove our main result.

\proclaim{\bf Theorem 3} For all\/ $n\in\Bbb Z^0$,\quad $(V_n,W_n)=(A_n,B_n)$. 
\endproclaim
\demo{\bf Proof} Since $(V_0,W_0)=(A_0,B_0)=(0,0)$, it suffices to show 
that for $n>0$, the numbers $V_n$, $W_n$ have the same inductive formation 
laws as the numbers $A_n$, $B_n$. By Lemma~2, $W_n=sV_n+tn$, the same as 
the formation rule for $B_n$ given in (2). It remains only to show that 
$V_n=\mex S$, where $S=\{V_i,W_i : i<n\}$. Suppose that $\mex S=W_j$. 
Clearly $j\ge n$. But then $V_j\in S$, since $V_j<W_j$. Hence $j<n$, a 
contradiction. 

Now $\bigcup_{i=1}^{\infty} V_i$ and $\bigcup_{i=1}^{\infty} W_i$ are 
complementary, since every positive integer has precisely one 
representation in $\Cal U$, either ending in an even number of $0$s, 
as the $V_n$ do, or in an odd number, as the $W_n$ do, since the $W_n$  
are a left shift of the $V_n$. Therefore if $\mex S\ne W_j$, we must 
have $\mex S=V_n$, so the formation laws are the same.\bqed\enddemo\bigskip

\centerline{{\bf 6. The End of the Games Story}}\medskip

Theorem~3 enables us to decide our main question, whether or not 
our heap games have a polynomial strategy. Given any game position 
$(x,y)$ with $0<x\le y$, compute $R(x)$ using the greedy algorithm 
mentioned in the first paragraph of \S4. If $R(x)$ ends in an odd number 
of $0$s, then $x=B_n$ for some $n>0$. Then move $y\rightarrow A_n$, 
where $R(B_n)=LR(A_n)$. If $R(x)$ ends in an even number 
of $0$s, then $x=A_n$ for some $n>0$. We can also test the relative 
size of $y$ and $B_n$, since $R(B_n)=LR(A_n)$. This information 
suffices for deciding the game, as indicated in the proof of Theorem~1 
(and used again in \S3). So the complexity of this computation 
is, up to a multiplicative constant, that of computing $R(x)$. 

The recurrence $u_n=ru_{n-1}+su_{n-2}$ has characteristic polynomial 
$x^2-rx-s=0$, with roots $\alpha=(r+\sqrt{r^2+4s})/2,\;$ 
$\beta=(r-\sqrt{r^2+4s})/2$. Since $s>0$, we have $\alpha>1$. 
Since $s=r-t+1\le r$, we have $0<-\beta\le(\sqrt{(r+2)^2-4}-r)/2<1$, 
so $|\beta|<1$. Therefore $u_n=E(c\alpha^n)$ for some constant $c>0$, 
where $E(v)$ is the nearest integer to the real number $v$. It 
follows that $n=O(\log x)$ bases of $\Cal U$ suffice for computing 
the strategy. Thus this strategy is in fact {\sl linear\/} in the 
input size.\bigskip

\centerline{{\bf 7. Yet Another Class of Sequences}}\medskip

In addition to the class of sequences $\{A_n\}$ and $\{B_n\}$ defined 
in (2), we now define another class of three sequences, $Q=\{Q_n\}$, 
$\{A'_n\}$, $\{B'_n\}$ $(n\in\Bbb Z^0)$, also depending on positive 
integer parameters $s$, $t$. 

\roster
\item"{(a)}" $Q_n=Q_m$ if $n=tQ_m+sm$ and $Q_m$ has already occurred 
precisely once; else 
$$Q_n=\mex\{Q_m:0\le m<n\}.\tag 6$$

\item"{(b)}" $A'_n=$ smallest $k$ such that $Q_k=n$. 

\item"{(c)}" $B'_n=$ largest $k$ such that $Q_k=n$. 
\endroster\medskip

Our main purpose here is to show that, despite the different definitions 
of the sequences, we actually have $A'_n=A_n$ and $B'_n=B_n$ for 
all $n\in\Bbb Z^0$.\medskip 

Partition $Q$ into subsequences $Q^1=\{Q^1_n\}$, 
$Q^2=\{Q^2_n\}$, where $Q^1$ consists of all the terms $Q_n=Q_m$ with 
smallest $m$, and $Q^2$ consists of the same terms, but with largest 
$n$.\medskip

{\bf Example.} For $s=2$, $t=1$, Table~3 lists the first few terms 
of these sequences.\medskip 

\midinsert
\topcaption{Table 3. The beginning terms of the five sequences for 
$s=2$, $t=1$.}
\endcaption
$$\vcenter{\halign{\tabskip 5pt\hfill #&\hfill #&\hfill #&\hfill #&
\hfill #&\hfill #&\hfill #&\hfill #&\hfill #&\hfill #&\hfill #&\hfill #&
\hfill #&\hfill #&\hfill #&\hfill #&\hfill #&\hfill #&\hfill #&\hfill #&
\hfill #&\hfill #&\hfill #\cr
$n$&\a&0&~1&~2&3&4&5&6&7&8&9&10&11&12&13&14&15&16&17&18&19&20\cr
\noalign{\vskip 3pt\hrule\vskip 3pt}
$Q_n$&\a&0&1&2&1&3&4&2&5&6&7&8&3&9&10&4&11&12&13&14&5&15\cr
$Q_n^1$&\a&0&1&2&&3&4&&5&6&7&8&&9&10&&11&12&13&14&&15\cr
$Q_n^2$&\a&0&&&1&&&2&&&&&3&&&4&&&&&5&\cr
$A^{\prime}_n$&\a&0&1&2&4&5&7&8&9&10&12&13&15&16&17&18&20&21&23&24&26&\cr
$B^{\prime}_n$&\a&0&3&6&11&14&19&22&25&&&&&&&&&&&&&\cr}}$$
\endinsert


It is convenient to precede the proof with two Lemmas. 

\proclaim{\bf Lemma 3} $(i)$ Let\/ $Q^2_i=r$, $Q^2_j=r+1$ be any two 
consecutive terms of\/ $Q^2$. Then\/ $j-i\ge 2$. 

\qquad\qquad $(ii)$ Let\/ $Q^1_i=r$, $Q^1_j=r+1$ be any two 
consecutive terms of\/ $Q^1$. Then\/ $j-i\le 2$.\endproclaim

\demo{\bf Proof} (i) We have $Q^2_i=r$ if $i=tQ_m+sm$ for some $m<i$, and 
$Q^2_j=r+1$ if $j=tQ_n+sn$ for some $m<n<i+1$, and $Q_m$, $Q_n$ occurred 
precisely once before. Then $j-i=t(Q_n-Q_m)+s(n-m)$. 
We clearly have $Q_n>Q_m$. Therefore $j-i\ge t+s\ge 2$. 

\qquad\quad (ii) The $\mex$ property (6) implies $j=i+1$, unless 
$i+1=tQ_m+sm$ for some $m<i+1$, where $Q_m$ appeared precisely once 
before. In this case $Q_{i+1}\in Q^2$. Part (i) implies that in this 
latter case $Q_{i+2}\in Q^1$, so then $j=i+2$.\bqed\enddemo

\proclaim{\bf Lemma 4} $A'_{r+1}-A'_r\in\{1,2\}$ for all\/ $r\in\Bbb Z^+$.
\endproclaim

\demo{\bf Proof} $A'_r$ is the smallest $i$ such that $Q_i=r$, and 
$A'_{r+1}$ is the smallest $j$ such that $Q_j=r+1$. The minimality of 
$i$ and $j$ means that $Q_i=Q^1_i$, $Q_j=Q^1_j$, and $Q^1_i=r$, 
$Q^1_j=r+1$ are consecutive. The result now follows from Lemma~3~(ii).
\bqed\enddemo

We are now ready to prove 

\proclaim{\bf Theorem 4} For every\/ $s, t\in\Bbb Z^+$ we have 
$A'_n=A_n$, $B'_n=B_n$ for all\/ $n\in\Bbb Z^+$, where\/ $A_n$, 
$B_n$ are defined in\/ $(2)$.\endproclaim

\demo{\bf Proof} Induction on $n$. By (2), $A_1=1$. Also $Q_1=1$ 
implies $A'_1=1$. Suppose we already showed that $A'_i=A_i$ for all 
$i\le n$. 

By Lemma~4, $A'_{n+1}-A'_n\in\{1,2\}$. In the proof 
of Theorem~2 (\S3), we showed that also $A_{n+1}-A_n\in\{1,2\}$ for 
all $n\in\Bbb Z^+$. 

Case (i).~~ $A_{n+1}=A_n+1$. Then by (2), $A_n+1=sA_m+tm$ for no 
$m\in\Bbb Z^+$. Suppose that $A'_{n+1}=A'_n+2$. Let $A'_n=k$. Then 
$A'_{n+1}=k+2$, $Q_k=n$, $Q_{k+2}=n+1$; and $Q_{k+1}=Q_{A'_n+1}$ has 
the property that $A'_n+1=A_n+1=tQ_m+sm$, where $Q_m$ has already 
occurred precisely once before. Thus if $Q_m=r$, then $A'_r=m$. 
Thus $A_n+1=tQ_{A'_n}+sA'_r=tr+sA'_r=tr+sA_r$ by the induction 
hypothesis, which is a contradiction. 

Case (ii).~~ $A_{n+1}=A_n+2$. The argument is similar to that of 
Case~(i), therefore it is omitted. Thus $A'_n=A_n$ for all $n\in\Bbb Z^+$. 

Now $B'_n=k$ if $Q_k=n$ and $Q_k$ occurred precisely once as some 
$Q_j$. It follows that for every $k\in\Bbb Z^+$ we have either 
$B'_n=k$ or $A'_n=k$, but not both. Hence $A'$, $B'$ are complementary 
sets, where $A'=\bigcup_{n=1}^{\infty} A'_n$, 
$B'=\bigcup_{n=1}^{\infty} B'_n$. The same was shown for 
$A=\bigcup_{n=1}^{\infty} A_n$, $B=\bigcup_{n=1}^{\infty} B_n$ at the 
beginning of the proof of Theorem~1 (\S2). It follows that also 
$B'_n=B_n$ for all $n\in\Bbb Z^+$.\bqed\enddemo

{\bf Notes.}
\roster
\item The special case $s=2$, $t=1$ of the second class (without 
$B'_n$) is listed in Neil Sloane's database of sequences \cite{Slo98}, 
sequence numbers 26366 $(Q_n)$, 26367 $(A'_n)$, ascribed to Clark 
Kimberling. We have not found a reference to other sequences of 
these families in \cite{Slo98}. In the definition of sequence 26366, 
``$a(n)=a(m)$ if  $m$ has already occurred exactly once$\dots$'', 
$m$ should presumably be replaced by $a(m)$.

\item We have $Q_0=0$, $Q_1=1$ for all $s,t\in\Bbb Z^0$; and 
$Q_2=1$ for $s=t=1\ $, $Q_2=2$ for all $s,t$ with $s+t>2$.  

\item The definition of the second class of sequences and the proof that 
both classes are identical, throws some light on the properties of 
both.
\endroster\bigskip 

\centerline{{\bf 8. Epilogue}}\medskip

The heap games proposed and analysed here belong to the family 
of {\it succinct\/} games, so named because their input size 
is succinct: $O(\log n)$ rather than $O(n)$. Often an extra 
effort is required for showing that such games are polynomial, 
i.e., have a polynomial strategy, because not more than $O(\log n)$ 
computation steps can be used. Different families of succinct 
games seem to require different methods of strategy computations. 

For example, in {\sl octal\/} games, invented by Guy and Smith 
\cite{GuSm56}, a linearly ordered string of beads may be split 
and or reduced according to rules encoded in octal. See also 
\cite{BCG82, \rm Ch.~4}, \cite{Con76, \rm Ch.~11}. The standard 
method for showing that an octal game is polynomial, is to demonstrate 
that its {\it Sprague-Grundy\/} function (the $0$s of which constitute 
the set of $P$-positions) is periodic. Periodicity has been 
established for a number of octal games. Some of the periods and or 
preperiods may be very large; see \cite{GaPl89}. Another way 
to establish polynomiality is to show that the Sprague-Grundy 
function values obey some other simple rule, such as forming 
an arithmetic sequence, as for Nim. 

For the present class of heap games, polynomiality was established 
by a nonstandard method. An arithmetic procedure, based on a class 
of special numeration systems, was the key to polynomiality. It 
appears that at this stage in the development of combinatorial game 
theory, there is no unified method for establishing polynomiality. 
But this malady seems to be common to most of discrete mathematics. 
Some might not even call it a malady, but consider it to be a feature 
inherent in the nature of mathematics. 

In \cite{YaYa67}, the special case of the Zeckendorf numeration system 
\cite{Zec72} was used to give one of the characterizations of the 
$P$-positions of Wythoff's game $(s=t=1)$. This method was extended in 
\cite{Fra82} for the generalized Wythoff game introduced there 
$(s=1$, $t\ge 1)$. In both cases, the bases of the 
numeration system were the numerators of the simple continued expansion 
of $\alpha$, were $\alpha$ is such that $A_n=\lfloor n\alpha\rfloor$ 
for all $n\ge 0$. The interesting aspect is that despite the fact that 
such $\alpha$ doesn't exist for $s>1$ (Theorem~2, \S3), the polynomial 
characterization based on special numeration systems nevertheless 
does exist. We also remark that it would be of interest to compute 
the Sprague-Grundy function for these heap games. For Wythoff's game 
this seems to be quite difficult, but this fact says nothing about the 
case $s>1$. 

In \cite{BoFr81} it is shown that a sequence $\{A_n\}$ is spectral 
(defined in the proof of Theorem~2), if and only if 
$|(A_{n+i}-A_n)-(A_{m+i}-A_m)|\le 1$ for all $i,m,n\ge 1$. Another 
motivation for the present paper was to extend this condition, namely 
to create and characterize sequences satisfying 
$|(A_{n+i}-A_n)-(A_{m+i}-A_m)|\le 2$. Vera S\'os told me that she has 
also been interested in this question. For the subfamily $s=t$ of the 
sequences $\{A_n\}$ defined in (2) (\S2), we have perhaps 
$|(A_{n+i}-A_n)-(A_{m+i}-A_m)|\le s$. And if this is true, does also 
the converse hold, namely, does $|(A_{n+i}-A_n)-(A_{m+i}-A_m)|\le s$ 
imply (2) with $s=t$? Investigation of the full family 
of these sequences (any $s,t\in\Bbb Z^+$), is of independent interest. 
In \S7 we defined a class of sequences, and demonstrated its equivalence 
with the class of sequences defined in (2). 

Not always is a succinct game more difficult than ``its nonsuccinct 
version''! We illustrate this with 
the game {\it vertex Kayles\/}. Given a finite (undirected) graph $G$. 
A move is to label an as yet unlabeled vertex not adjacent to any 
labeled vertex. The player first unable to play loses, and the opponent 
wins. A partizan variation is called {\it bigraph vertex Kayles\/}. 
Both versions have been proved Pspace-hard in \cite{Sch78}. 
If $G$ is a path, the resulting succinct game, known as {\it Kayles\/}, 
is actually polynomial! It is the octal game $0.137$ --- see
\cite{GuSm56}, \cite{BCG82}, \cite{Con76}. Incidentally, there is a 
large ``no-man's-land'' of games lying in between the polynomial $0.137$ 
and the Pspace-hard vertex Kayles, and it would be of interest to reduce
the boundary area. 
 
Finally, the family of combinatorial games consists, roughly, of 
two-player games with perfect information (no hidden information 
as in some card games), no chance moves (no dice) and outcome 
restricted to (lose, win). These games are {\sl completely determined\/}, 
so one of their main mathematical interests is in bounding the complexity 
of their strategies. This explains why we talked so much about efficiency 
of strategy computation in this paper.\bigskip

\centerline{{\bf References}}\medskip

1.\ [BCG82] E.R.\ Berlekamp, J.H.\ Conway and R.K.\ Guy \(1982\), {\it Winning
Ways\/} (two volumes), Academic Press, London.

2.\ [BoFr81] M.\ Boshernitzan and A.S.\ Fraenkel \(1981\), 
Nonhomogeneous spectra of numbers,  
{\it Discr.\ Math.\/} {\bf 34}, 325--327

3.\ [BoFr84] M.\ Boshernitzan and A.S.\ Fraenkel \(1984\), 
A linear algorithm for nonhomogeneous spectra of numbers, 
{\it J.\ of Algorithms\/} {\bf 5}, 187--198.

4.\ [Con76] J.H.\ Conway \(1976\), {\it On Numbers and Games}, Academic Press,
London.

5.\ [Cox53] H.S.M.\ Coxeter \(1953\), The golden section, phyllotaxis 
and Wythoff's game, {\it Scripta Math.\/} {\bf 19}, 135--143.

6.\ [Fra82] A.S.\ Fraenkel \(1982\), How to beat your Wythoff games' 
opponents on three fronts, {\it Amer.\ Math.\ Monthly\/} {\bf 89}, 353--361.

7.\ [Fra85] A.S.\ Fraenkel \(1985\), Systems of numeration,
{\it Amer.\ Math.\ Monthly\/} {\bf 92}, 105--114.

8.\ [GaPl89] A.\ Gangolli and T.\ Plambeck \(1989\), A note on periodicity 
in some octal games, {\it Internat.\ J.\ Game Theory\/} {\bf 18}, 311--320.

9.\ [GLL78] R.L.\ Graham, S.\ Lin and C.-S.\ Lin \(1978\), Spectra of numbers, 
{\it Math.\ Mag.\/} {\bf 51}, 174--176.

10. [GuSm56] R.K.\ Guy and C.A.B.\ Smith \(1956\), The $G$-values of 
various games, {\it Proc.\ Camb.\ Phil.\ Soc.\/} {\bf 52}, 514--526.

11. [Sch78] T.J.\ Schaefer \(1978\), On the complexity of some two-person
perfect-in\-for\-ma\-tion games, {\it J.\ Comput.\ System Sci.} {\bf 16},
185--225.

12. [Slo98] N.J.A.\ Sloane \(1998\), Sloane's On-Line Encyclopedia of 
Integer Sequences, http://www.research.att.com/\~\! njas/sequences/ . 

13. [Wyt07] W.A.\ Wythoff \(1907\), A modification of the game of Nim, 
{\it Nieuw Arch.\ Wisk.\/} {\bf 7}, 199--202.

14. [YaYa67] A.M.\ Yaglom and I.M.\ Yaglom \(1967\), {\it Challenging 
Mathematical Problems with Elementary Solutions\/}, translated by 
J. McCawley, Jr., revised and edited by B. Gordon, Vol.~II, Holden-Day, 
San Francisco.

15. [Zec72] E.\ Zeckendorf \(1972\), Repr\'esentation des nombres naturels 
par une somme de nombres de Fibonacci ou de nombres de Lucas, 
{\it Bull.\ Soc.\ Roy.\ Sci.\ Li\`ege\/} {\bf 41}, 179--182.

\enddocument